\documentclass[english,11pt]{article}
\usepackage{amssymb}
\usepackage{amsmath,amsfonts,amssymb,amsthm}
\usepackage{fullpage}
\usepackage{latexsym}
\usepackage{xspace}
\usepackage{babel}
\usepackage{paralist}
\usepackage{algorithm}
\usepackage{algorithmic}

\newtheorem{theorem}{Theorem}
\newtheorem{lemma}{Lemma}

\begin{document}

\title{On Rainbow-$k$-Connectivity of Random Graphs\footnotemark[1]}
\author{Jing He\footnotemark[2] \and Hongyu Liang\footnotemark[3]}

\renewcommand{\thefootnote}{\fnsymbol{footnote}}

\footnotetext[1]{This work was supported in part by the National
Basic Research Program of China Grant 2011CBA00300, 2011CBA00301,
and the National Natural Science Foundation of China Grant 61033001,
61061130540, 61073174.}

\footnotetext[2]{nstitute for Interdisciplinary Information
Sciences, Tsinghua University, e-mail: he-j08@mails.tsinghua.edu.cn}

\footnotetext[3]{nstitute for Interdisciplinary Information
Sciences, Tsinghua University, e-mail:
lianghy08@mails.tsinghua.edu.cn}

\date{December 15, 2010}
\maketitle

\begin{abstract}
A path in an edge-colored graph is called a \emph{rainbow path} if
the edges on it have distinct colors. For $k\geq 1$, the
\emph{rainbow-$k$-connectivity} of a graph $G$, denoted $rc_k(G)$,
is the minimum number of colors required to color the edges of $G$
in such a way that every two distinct vertices are connected by at
least $k$ internally vertex-disjoint rainbow paths. In this paper,
we study rainbow-$k$-connectivity in the setting of random graphs.
We show that for every fixed integer $d\geq 2$ and every $k\leq
O(\log n)$, $p=(\log n)^{1/d}/n^{(d-1)/d}$ is a sharp threshold
function for the property $rc_k(G(n,p))\leq d$. This substantially
generalizes a result in [Y. Caro, A. Lev, Y. Roditty, Z. Tuza, and
R. Yuster, On rainbow connection, Electron. J. Comb., 15, 2008],
stating that $p=\sqrt{\log n/n}$ is a sharp threshold function for
the property $rc_1(G(n,p))\leq 2$. As a by-product, we obtain a
polynomial-time algorithm that makes $G(n,p)$ rainbow-$k$-connected
using at most one more than the optimal number of colors with
probability $1-o(1)$, for all $k\leq O(\log n)$ and
$p=n^{-\epsilon(1\pm o(1))}$ for any constant $\epsilon\in[0,1)$.
\end{abstract}

\section{Introduction}\label{sec:intro}
All graphs considered in this paper are finite, simple, undirected
and contain at least 2 vertices. We follow the notation and
terminology of \cite{b98}. The following notion of
\emph{rainbow-$k$-connectivity} was proposed by Chartrand et al.
\cite{cjmz08,cjmz09} as a strengthening of the canonical
connectivity concept in graphs. Given an edge-colored graph $G$, a
path in $G$ is called a \emph{rainbow path} if its edges have
distinct colors. For an integer $k\geq 1$, an edge-colored graph is
called \emph{rainbow-$k$-connected} if any two different vertices of
$G$ are connected by at least $k$ internally vertex-disjoint rainbow
paths. The \emph{rainbow-$k$-connectivity} of $G$, denoted by
$rc_k(G)$, is the minimum number of colors required to color the
edges of $G$ to make it rainbow-$k$-connected. Note that such
coloring does not exist if $G$ is not $k$-vertex-connected, in which
case we simply let $rc_k(G)=\infty$. When $k=1$ it is alternatively
called \emph{rainbow-connectivity} or \emph{rainbow connection
number} in literature, and is conventionally written as $rc(G)$ with
the subscript $k$ dropped.

Besides its theoretical interest as being a natural combinatorial
concept, rainbow connectivity also finds applications in networking
and secure message transmitting \cite{cfmy09,e07,s10}. The following
motivation is given in \cite{cfmy09}: Suppose we want to route
messages in a cellular network such that each link on the route
between two vertices is assigned with a distinct channel. Then the
minimum number of used channels is exactly the rainbow-connectivity
of the underlying graph.

Some easy observations regarding rainbow-$k$-connectivity include
that $rc_k(G) = 1$ if and only if $k=1$ and $G$ is a clique, that
$rc(G)\leq n-1$ for all connected $G$, and that $rc(G)=n-1$ if and
only if $G$ is a tree, where $n$ is the number of vertices in $G$.
Chartrand et al. \cite{cjmz08} determined the rainbow-connectivity
of several special classes of graphs, including complete
multipartite graphs. In \cite{cjmz09} they investigated
rainbow-$k$-connectivity in complete graphs and regular complete
bipartite graphs. The extremal graph-theoretic aspect of
rainbow-connectivity was studied by Caro et al. \cite{clrty08}, who
proved that $rc(G)=O_{\delta}(n\log \delta/\delta)$ with $\delta$
being the minimum degree of $G$. This tradeoff was later improved to
$rc(G)<20n/\delta$ by Krivelevich and Yuster~\cite{ky10}, and was
recently shown to be $rc(G)\leq 3n/(\delta+1)+3$ by Chandran et
al.~\cite{cdrv10} which is essentially tight. Chakraborty et al.
\cite{cfmy09} studied the computational complexity perspective of
this notion, proving among other results that given a graph $G$
deciding whether $rc(G)=2$ is NP-complete.

Another important setting that has been extensively explored for
studying various graph concepts is the Erd\H os-R\'enyi random graph
model $G(n,p)$~\cite{er60}, in which each of the ${n\choose 2}$
pairs of vertices appears as an edge with probability $p$
independently from other pairs. we say an event $\mathcal{E}$
happens \emph{almost surely} if the probability that it happens
approaches 1 as $n \rightarrow \infty$, i.e.,
$\textbf{Pr}[\mathcal{E}]=1-o_n(1)$. We will always assume that $n$
is the variable that tends to infinity, and thus omit the subscript
$n$ from the asymptotic notations. For a graph property $P$, a
function $p(n)$ is called a \emph{threshold function} of $P$ if:
\begin{itemize}
\item for every $r(n)=\omega(p(n))$, $G(n,r(n))$ almost surely satisfies $P$; and
\item for every $r'(n)=o(p(n))$, $G(n,r'(n))$ almost surely does not satisfy $P$.
\end{itemize}
Furthermore, $p(n)$ is called a \emph{sharp threshold function} of
$P$ if there exist two positive constants $c$ and $C$ such that:
\begin{itemize}
\item for every $r(n)\geq C\cdot p(n)$, $G(n,r(n))$ almost surely satisfies $P$; and
\item for every $r'(n)\leq c\cdot p(n)$, $G(n,r'(n))$ almost surely does not satisfy $P$.
\end{itemize}
Clearly a sharp threshold function of a graph property is also a
threshold function of it; yet the converse may not hold, e.g., the
property of containing a triangle \cite{b85}.

It is known that every non-trivial monotone graph property possesses
a threshold function \cite{bt86,ek96}. Obviously for every $k,d$,
the property $rc_k(G)\leq d$ is monotone, and thus has a threshold.
Caro et al. \cite{clrty08} proved that $p=\sqrt{\log n/n}$ is a
sharp threshold function for the property $rc_1(G(n,p))\leq 2$.
In this paper, we significantly extend their result by establishing
sharp thresholds for the property $rc_k(G(n,p))\leq d$ for all
constants $d$ and logarithmically increasing $k$. Our main theorem
is as follows.

\begin{theorem}\label{thm:thr}
Let $d\geq 2$ be a fixed integer and $k=k(n)\leq O(\log n)$. Then
$p=(\log n)^{1/d}/n^{(d-1)/d}$ is a sharp threshold function for the
property $rc_k(G(n,p))\leq d$.
\end{theorem}

We also investigate rainbow-$k$-connectivity from the algorithmic
point of view. The NP-hardness of determining $rc(G)$ is shown by
Chakraborty et al. \cite{cfmy09}. We show that the problem (even the
search version) becomes easy in random graphs, by designing an
algorithm for coloring random graphs to make it
rainbow-$k$-connected with near-optimal number of colors.

\begin{theorem}\label{thm:alg}
For any constant $\epsilon\in [0,1)$, $p=n^{-\epsilon(1\pm o(1))}$
and $k\leq O(\log n)$, there is a randomized polynomial-time
algorithm that, with probability $1-o(1)$, makes $G(n,p)$
rainbow-$k$-connected using at most one more than the optimal number
of colors, where the probability is taken over both the randomness
of $G(n,p)$ and that of the algorithm.
\end{theorem}

Our result is quite strong, since almost all natural edge
probability functions $p$ encountered in various scenarios satisfy
$p=n^{-\epsilon(1\pm o(1))}$ for some $\epsilon>0$. Note that
$G(n,n^{-\epsilon})$ is almost surely disconnected when
$\epsilon>1$~\cite{er60}, which makes the problem become trivial. We
therefore ignore these cases.

In Section~\ref{sec:rc_random} we present the proof of
Theorem~\ref{thm:thr}, and in Section~\ref{sec:find_rc} we show the
correctness of Theorem~\ref{thm:alg}.

\section{Threshold of Rainbow-$k$-Connectivity}\label{sec:rc_random} This section is devoted to proving Theorem~\ref{thm:thr}. Throughout the paper ``$\ln$'' denotes the natural logarithm, and ``$\log$'' denotes the logarithm to the base 2. Hereafter we assume $d\geq 2$ is a
fixed integer, $c_0\geq 1$ a positive constant, and $k=k(n)\leq
c_0\log n$ for all sufficiently large $n$. To establish a sharp
threshold function for a graph property the proof should be
two-fold. We first show the easy direction.

\begin{theorem}\label{thm:lb}
$rc_k(G(n,(\ln n)^{1/d}/n^{(d-1)/d}))\geq d+1$ almost surely holds.
\end{theorem}

We need the following fact proved by Bollob\'as \cite{b81}.

\begin{lemma}[Restatement of part of Theorem 6 in~\cite{b81}] Let $c$ be a
positive constant and $d\geq 2$ a fixed integer. Let
$p'=(\ln(n^2/c))^{1/d}/n^{(d-1)/d}$. Then,
$$ \lim_{n\rightarrow \infty}\emph{\textbf{Pr}}[G(n,p')~\textrm{has diameter at most~} d]=e^{-c/2}.$$
\end{lemma}

\begin{proof}[of Theorem~\ref{thm:lb}]
Fix an arbitrary $\epsilon>0$ and choose a constant $c>0$ so that
$e^{-c/2}<\epsilon/2$. Let $p'=(\ln(n^2/c))^{1/d}/n^{(d-1)/d}$ and
$p=(\ln n)^{1/d}/n^{(d-1)/d}$. Clearly $p\leq p'$ for all $n>c$.

By Lemma 1 and the definition of limits, there exists an $N_1>0$
such that for all $n>N_1$, $\textbf{Pr}[G(n,p')~\textrm{has diameter
at most~} d]<e^{-c/2}+\epsilon/2<\epsilon,$ by our choice of $c$.
Thus, for every $n>\max\{c,N_1\}$,
\begin{eqnarray*}
\textbf{Pr}[G(n,p)~\textrm{has diameter at most~} d]\leq
\textbf{Pr}[G(n,p')~\textrm{has diameter at most~} d]<\epsilon.
\end{eqnarray*}

Due to the arbitrariness of $\epsilon$, this implies that the
probability of $G(n,p)$ having diameter at most $d$ is $o(1)$. This
completes the proof of Theorem~\ref{thm:lb}, since the
rainbow-$k$-connectivity of a graph is at least as large as its
diameter.
\end{proof}

We are left with the other direction stated below. Fix
$C=2^{20}\cdot c_0$.

\begin{theorem}\label{thm:ub}
$rc_k(G(n,C(\log n)^{1/d}/n^{(d-1)/d}))\leq d$ almost surely holds.
\end{theorem}

The key component of our proof of Theorem~\ref{thm:ub} is the
following theorem.

\begin{theorem}\label{thm:disjoint_path}
With probability at least $1-n^{-\Omega(1)}$, every two different
vertices of $G(n,C(\log n)^{1/d}/n^{(d-1)/d})$ are connected by at
least $2^{10d}c_0\log n$ internally vertex-disjoint paths of length
exactly $d$.
\end{theorem}

Before demonstrating Theorem~\ref{thm:disjoint_path}, we show how
Theorem~\ref{thm:ub} follows from it.

\begin{proof}[of Theorem~\ref{thm:ub}]
Let $G$ be an instance of $G(n,C(\log n)^{1/d}/n^{(d-1)/d})$ for
which the condition in Theorem~\ref{thm:disjoint_path} holds; that
is, every two different vertices of $G$ have at least $c_1\log n$
internally vertex-disjoint paths of length $d$ connecting them,
where $c_1:=2^{10d}c_0$. To establish Theorem~\ref{thm:ub} it
suffices to prove that $rc_k(G)\leq d$ for every such $G$, since by
Theorem~\ref{thm:disjoint_path} the condition holds with probability
at least $1-n^{-\Omega(1)}=1-o(1)$.

Let $S=\{1,2,\ldots,d\}$ be a set of $d$ distinct colors. We
randomly color the edges of $G$ with colors from $S$. Fix two
distinct vertices $u$ and $v$. Let $P$ be a path of length $d$
connecting $u$ and $v$. The probability that $P$ becomes a rainbow
path under the random coloring is
$$q:=d!/d^d\geq (d/e)^d/d^d\geq 4^{-d},$$
by Stirling's formula. Since there are at least $c_1\log n$ such
paths and they are all edge-disjoint, we can upper-bound the
probability that at most $k-1$ of them become rainbow paths by
\begin{eqnarray*}
& &{c_1\log n \choose k-1}(1-q)^{c_1\log n-(k-1)}\\
&\leq& {c_1\log n\choose c_0\log n}(1-4^{-d})^{(c_1-c_0)\log n}\\
&\leq& 2^{c_1\log n\cdot H(c_0/c_1)}\cdot
2^{-4^{-d}(c_1-c_0)\log n}\\
&=& n^{-(4^{-d}(c_1-c_0)-c_1\cdot H(c_0/c_1))},
\end{eqnarray*}
where the second inequality follows from the fact that
$${m\choose
\alpha m}\leq 2^{m\cdot H(\alpha)}$$ for all constants
$\alpha\in(0,1)$ and sufficiently large $m$, $H$ being the binary
entropy function defined as
$$H(\epsilon)=\epsilon\log(1/\epsilon)+(1-\epsilon)\log(1/(1-\epsilon)),$$
and that $$1-x\leq e^{-x}\leq 2^{-x},~\textrm{for all~}x\geq 0.$$

It is easy to verify that $\log x\leq \sqrt{x}$ whenever $x\geq
100$. Also, since $1+x\leq e^{x}\leq 2^{2x}$, we have $\log(1+x)\leq
2x$ for all $x>-1$. Recalling that $c_1=2^{10d}c_0>200c_0$, we get
\begin{eqnarray*}
H(c_0/c_1)&=&(c_0/c_1)\log
(c_1/c_0)+(1-c_0/c_1)\log (1+c_0/(c_1-c_0))\\
&\leq& (c_0/c_1)\sqrt{c_1/c_0}+(1-c_0/c_1)\cdot
2c_0/(c_1-c_0) \\
&=&\sqrt{c_0/c_1}+2c_0/c_1 \leq 3\sqrt{c_0/c_1}\,.
\end{eqnarray*}
We thus have
\begin{eqnarray*}
4^{-d}(c_1-c_0)-c_1\cdot H(c_0/c_1) &\geq& 4^{-d}(c_1-c_0)- 3\sqrt{c_1c_0}\\
&=& 4^{-d}c_1(1-2^{-10d})- 3\sqrt{2^{-10d}\cdot c_1^2 }\\
&\geq& 2^{-2d-1}c_1 - 2^{-5d+2}c_1\\
&\geq& c_1\cdot 2^{-2d-2}\\
&=&c_0\cdot 2^{10d}\cdot 2^{-2d-2}>100\, ,
\end{eqnarray*}
since $c_0\geq 1$ and $d\geq 2$. Therefore, the probability that
there exist at most $k-1$ rainbow paths between $u$ and $v$ is at
most
$$n^{-(4^{-d}(c_1-c_0)-c_1\cdot H(c_0/c_1))}<n^{-100}.$$
By the Union Bound, with probability at least
$$1-{n\choose 2}n^{-100}\geq 1-n^{-90},$$
every two distinct vertices of $G$ have at least $k$ internally
vertex-disjoint rainbow paths connecting them. In particular, there
exists a $d$-coloring of the edges of $G$ under which $G$ becomes
$k$-rainbow-connected, implying that $rc_k(G)\leq d$. This concludes
the proof of Theorem~\ref{thm:ub}.
\end{proof}

We now prove Theorem~\ref{thm:disjoint_path}.

\begin{proof}[of Theorem~\ref{thm:disjoint_path}]
Let $p=C(\log n)^{1/d}/n^{(d-1)/d}$ where $C=2^{20}c_0$. Let $V$ be
the set of all vertices in $G(n,p)$. For every $S\subseteq V$ and
$u\in S$, let $\mathcal{A}(S,u)$ be the event that $u$ is adjacent
to at least $pn/10(=Cn^{1/d}(\log n)^{1/d}/10)$ distinct vertices in
$V\setminus S$. The following lemma is needed for our proof.
\begin{lemma}\label{lem:manynei}
For every $S,u$ such that $u\in S$ and $|S|\leq d\cdot
(pn/10)^{d-1}$,
$$\emph{\textbf{Pr}}[\mathcal{A}(S,u)]\geq 1-2^{-\Omega(n^{1/d})}.$$
\end{lemma}
\begin{proof}
Fix $S\subseteq V$ with $|S|\leq d\cdot (pn/10)^{d-1}$, and $u\in
S$. We have
$$|V\setminus S|\geq n-d\cdot(pn/10)^{d-1}=n-d\cdot (C/10)^{d-1}n^{(d-1)/d}(\log
n)^{(d-1)/d}\geq n/2,$$ for all sufficiently large $n$. Take $T$ to
be any subset of $V\setminus S$ of cardinality $n/2$. Let $X$ be the
random variable counting the number of neighbors of $u$ inside $T$.
It is obvious that $X$ can be expressed as the sum of $n/2$
independent random variables, each of which taking 1 with
probability $p$ and 0 with probability $1-p$. Thus
$\textbf{E}[X]=pn/2$. By the Chernoff-Hoeffding Bound (see e.g.
Theorem 4.2 of \cite{mr}), we have
\begin{eqnarray*}
\textbf{Pr}[X<(1-4/5)pn/2]\leq \exp(-(1/2)(4/5)^2
(pn/2))=2^{-\Omega(n^{1/d})},
\end{eqnarray*}
which gives precisely what we want.
\end{proof}

We now continue the proof of Theorem~\ref{thm:disjoint_path}. Fix
$u,v\in V, u\neq v$. Let $S_0=\{u\}$. A \emph{$t$-ary tree} with a
designated root is a tree whose non-leaf vertices all have exactly
$t$ children. Consider the following process of ``choosing'' a
$(pn/10)$-ary tree of depth $d-1$ rooted at $u$:
\begin{enumerate}
\item Let $i\leftarrow 1$ and $S_i\leftarrow \emptyset$. \item For
every vertex $w\in S_{i-1}$ (in an arbitrary order), choose $pn/10$
distinct neighbors of $w$ from the set $V\setminus(\{v\}\cup
\bigcup_{j=0}^{i}S_j)$, and add them to $S_i$. (Note that $S_i$ is
updated every time after the processing of a vertex $w$, so that one
vertex cannot be chosen and added to $S_i$ more than once. This
ensures that at the end of this step, $|S_i|=(pn/10)^{i}$.) \item
Let $i\leftarrow i+1$. If $i\leq d-1$ then go to Step 2, otherwise
stop.
\end{enumerate}

Of course the process may fail during Step 2, since with nonzero
probability $w$ will have no neighbor in $V\setminus(\{v\}\cup
\bigcup_{j=0}^{i}S_j)$. (In fact, with nonzero probability the graph
becomes empty.) However, noting that at any time during the process,
$$|\{v\}\cup \bigcup_{j=0}^{i}S_j|\leq
1+\sum_{j=0}^{d-1}(pn/10)^j\leq d\cdot (pn/10)^{d-1},\textrm{~for
all sufficiently large~}n,$$ we can deduce from
Lemma~\ref{lem:manynei} that every execution of Step 2 fails with
probability at most $2^{-\Omega(n^{1/d})}.$ Since Step 2 can be
conducted for at most $d\cdot (pn/10)^{d-1}$ times, we obtain that,
with probability at least
$$1-d\cdot (pn/10)^{d-1}\cdot 2^{-\Omega(n^{1/d})}=1-2^{-\Omega(n^{1/d})},$$
the process will successfully terminate. At the end of the process,
the sets $S_0,S_1,\ldots,S_{d-1}$ naturally induces a $(pn/10)$-ary
tree $T$ of depth $d-1$ rooted at $u$, with $S_i$ being the
collection of vertices in the $i$-th level. The number of leaves in
$T$ is exactly $|S_{d-1}|=(pn/10)^{d-1}$.

Now we assume that $T$ has been successfully constructed. Let $Y$ be
a random variable denoting the number of neighbors of $v$ inside
$S_{d-1}$. (Recall that $v\not\in S_{d-1}$.) It is clear that
$$\textbf{E}[Y]=p\cdot|S_{d-1}|=p^{d}n^{d-1}/10^{d-1}=10\cdot (C/10)^{d}\log n.$$
As before, using the Chernoff-Hoeffding Bound, we have
$$\textbf{Pr}[Y<(C/10)^{d}\log n]\leq \exp(-(1/2)(9/10)^2(C/10)^d\cdot 10\log n)\leq n^{-10}, $$
for our choice of $C$.

It is clear that each neighbor $v'$ of $v$ inside $S_{d-1}$ induces
a length-$d$ path between $u$ and $v$ (by simply combining the path
from $u$ to $v'$ in tree $T$ and the edge $(v',v)$). The problem is
that these paths may not be internally vertex-disjoint. We next
address this issue.

For every $w\in S_1$,  denote by $T_{w}$ the subtree of $T$ of depth
$d-2$ rooted at $w$. Call these $T_{w}$ \emph{vice-trees}. Clearly
every vice-tree contains $(pn/10)^{d-2}$ leaves.

The reason for defining such vice-trees is that, by simple
observations, any two leaves of $T$ that belong to different
vice-trees must correspond to edge-disjoint root-to-leaf paths in
$T$. Thus, to establish a large number of internally vertex-disjoint
paths between $u$ and $v$, it suffices to show that we can find many
neighbors of $v$ inside $S_{d-1}$ that belong to distinct
vice-trees.

For each vice-tree $T_{w}$, let $\mathcal{B}_w$ be the event that
$v$ has at least $10d$ neighbors inside the set of leaves of
$T_{w}$. Noting that $T_w$ has exactly $(pn/10)^{d-2}$ leaves, we
have

\begin{eqnarray*}
\textbf{Pr}[\mathcal{B}_w]&\leq& {(pn/10)^{d-2}\choose
10d}p^{10d}\\
&=&{(Cn^{1/d}(\log n)^{1/d}/10)^{d-2}\choose 10d}\cdot
\left(\frac{C(\log n)^{1/d}}{n^{(d-1)/d}}\right)^{10d}\\
&\leq&\left((Cn^{1/d}(\log n)^{1/d}/10)^{d-2}\cdot
\frac{C(\log n)^{1/d}}{n^{(d-1)/d}}\right)^{10d}\\
&=&C^{10d}(C/10)^{10d(d-2)}(\log n)^{10(d-1)}n^{-10}\\
&\leq &O(n^{-9}).
\end{eqnarray*}

By applying the Union Bound, we obtain
$$\textbf{Pr}[\bigvee_{w\in S_1}\mathcal{B}_w]\leq
(pn/10)\cdot O(n^{-9})\leq O(n^{-7}).$$

Combined with previous results, we deduce that with probability at
least
$$1-2^{-\Omega(n^{1/d})}-n^{-10}-O(n^{-7})\geq 1-O(n^{-6}),$$
the following three events simultaneously happen:
\begin{enumerate}
\item The tree $T$ is successfully constructed.

\item $v$ has at least $(C/10)^{d}\log n$ neighbors that are
leaves of $T$.

\item Every vice-tree $T_{w}$ contains at most $10d$ leaves that
are neighbors of $v$.
\end{enumerate}

When all these three events happen, we can choose
$((C/10)^{d}/(10d))\log n$ neighbors of $v$, every two of which come
from different vice-trees. This immediately leads to
$((C/10)^{d}/(10d))\log n$ length-$d$ internally vertex-disjoint
paths between $u$ and $v$, where, by our choice of $C=2^{20}c_0$,
$$((C/10)^{d}/(10d))\log n\geq 2^{10d}c_0\log n. $$

Using the Union Bound again gives us that, with probability at least
$$1-{n\choose 2}\cdot O(n^{-6})=1-n^{-\Omega(1)},$$
every two distinct vertices have at least $2^{10d}c_0\log n$
internally vertex-disjoint paths of length $d$ connecting them. The
proof of Theorem~\ref{thm:disjoint_path} is thus completed.

\end{proof}

\section{Rainbow-coloring Random Graphs}\label{sec:find_rc}
In this section we prove Theorem~\ref{thm:alg}.

\begin{proof}[of Theorem~\ref{thm:alg}]
First note that for every $G$ with at least 2 vertices, $rc_k(G)=1$
if and only if $k=1$ and $G$ is a clique, which can be easily
checked. Thus, in the following we assume w.l.o.g. that
$rc_k(G(n,p))\geq 2$.

It is easy to see that there exists a (unique) integer $d\geq 2$
such that $(d-2)/(d-1)\leq \epsilon < (d-1)/d$. We have
$p=\omega\left((\log n)^{1/d}/n^{(d-1)/d}\right)$, which, by
Theorem~\ref{thm:ub}, implies that $rc_k(G(n,p))\leq d$ almost
surely holds. Moreover, a scrutiny into the proof of
Theorem~\ref{thm:disjoint_path} tells us that for such $p$, a random
coloring of $G(n,p)$ using $d$ colors will make it
rainbow-$k$-connected with probability $1-n^{-\Omega(1)}$. Thus, it
suffices for us to show that with probability $1-o(1)$,
$rc_k(G(n,p))\geq d-1$. This is trivial for $d\leq 3$, since we have
assumed that $rc_k(G(n,p))\geq 2$. When $d\geq 4$, we have
$p=o\left((\log n)^{1/(d-2)}/n^{(d-3)/(d-2)}\right)$. Due to
Theorem~\ref{thm:lb}, $G(n,p)$ with such $p$ almost surely satisfies
$rc_k(G(n,p))\geq d-1$.

Hence, we have shown that with probability $1-o(1)$, a random
coloring with $d$ colors will make $G(n,p)$ rainbow-$k$-connected
and the number of colors used is at most one more than the optimum,
where the probability is taken over both the randomness of $G(n,p)$
and that of the algorithm. This completes the whole proof.
\end{proof}

\bibliographystyle{plain}
\bibliography{rc_threshold_full}

\end{document}